\documentclass[12pt]{article}
\usepackage[cp1251]{inputenc}
\usepackage{amsthm,amsfonts,amsmath,amscd,amssymb}
\usepackage{latexsym}
\usepackage[english]{babel}
\usepackage[T2A]{fontenc}
\usepackage{graphicx}
\usepackage{caption}
\usepackage{subcaption}
\usepackage{cite}

\begin{document}

\begin{large}

\centerline{\textbf{On Circular Tractrices in $\mathbb R^3$}}

\end{large}

\bigskip

\centerline{V.~Gorkavyy\footnote{e-mail: gorkaviy@ilt.kharkov.ua, vasylgorkavyy@gmail.com}  }

\centerline{(B.~Verkin  Institute for Low Temperature Physics and Engineering, Kharkiv, Ukraine)}

\centerline{A.~Sirosh\footnote{e-mail: alina.sirosh98@gmail.com} }

\centerline{(V.N.~Karazin Kharkiv National University, Kharkiv, Ukraine)}


{\small {\bf Abstract.} We explore geometric properties of circular analogues of tractrices and pseudospheres in $\mathbb R^3$.}



\bigskip

	\smallskip
	\section{Introduction}
	
	In 2000 Yuriy Aminov and Antony Sym settled the question whether one can extend the classical theory of Bianchi-Backlund transformations of pseudo-spherical surfaces in $\mathbb R^3$ to the case of pseudo-spherical surfaces in $\mathbb R^4$, see \cite{AminovSym}. This question, as well as its multi-dimensional generalizations, was addressed in a series of research papers \cite{G0} - \cite{G5}, where it was shown to be rather non-trivial and actually it still remains widely open.
	
	While studying the problem we (re)discovered a particular family of spatial curves in $\mathbb R^3$ called {\it circular tractrices}, which can be used for constructing novel examples of pseudo-spherical surfaces in $\mathbb R^4$ similar to the Beltrami and Dini surfaces in $\mathbb R^3$ \cite{GStepanova}. Despite possible applications to the theory of pseudo-spherical surfaces, the circular tractrices themselves turn out to be of independent interest. The aim of this research note is to survey beautiful geometric properties of circular tractrices with particular emphasis on justifying the use of the terms {\it tractrix} and {\it circular}.
	
	Let us introduce the principal hero of our story.
	
	Consider the three-dimensional Euclidean space $\mathbb R^3$ endowed with Cartesian coordinates $(x^1,x^2,x^3)$.
	
{\bf Definition.} A {\it circular tractrix}  is a curve in $\mathbb R^3$ represented by
		\begin{eqnarray}\label{pv}\left\{
			\begin{array}{l}
				x^1 = \xi_1 \cos\frac{t}{R}\, +\, \xi_2 \sin\frac{t}{R},\\
				x^2 = - \xi_2 \cos\frac{t}{R}\, +\, \xi_1 \sin\frac{t}{R},\\
				x^3 = \xi_3,
			\end{array}\right. \quad\quad\quad  t\in\mathbb R,
		\end{eqnarray}
		where $R>0$ is a fixed constant, and $\xi_1(t)$, $\xi_2(t)$, $\xi_3(t)$ are functions given explicitly by the following formulae depending on whether $R$ is greater, equal or less than 1 respectively:
		\begin{eqnarray}\label{large}R>1)\quad \xi_1 = \frac{(R-\frac{1}{R})\cosh\lambda t}{\frac{c_1}{R}+\cosh\lambda t}, \quad \xi_2 = \frac{\lambda\sinh\lambda t}{\frac{c_1}{R}+\cosh\lambda t}, \quad \xi_3 =\frac{\lambda c_2}{\frac{c_1}{R}+\cosh\lambda t},
		\end{eqnarray}
		where $\lambda=\frac{\sqrt{R^2-1}}{R}$, and $c_1$, $c_2$ are arbitrary constants subject to $c_1^2+c_2^2=1$;
		\begin{eqnarray}\label{unit}R=1)\quad \quad \xi_1 = \frac{2}{c_1+t^2}, \quad \xi_2 =  \frac{2t}{c_1+t^2}, \quad \xi_3 = \frac{c_2}{c_1+t^2},
		\end{eqnarray}
		where $c_1$, $c_2$ are arbitrary constants subject to $4(c_1-1) = c_2^2$;
		\begin{eqnarray}\label{small}R<1)\quad
			\xi_1=\frac{(R-\frac{1}{R})\cos\lambda t}{\frac{c_1}{R}+\cos\lambda t}, \quad \xi_2=-\frac{\lambda\sin\lambda t}{\frac{c_1}{R}+\cos\lambda t}, \quad \xi_3=\frac{\lambda c_2}{\frac{c_1}{R}+\cos\lambda t},
		\end{eqnarray}
		where  $\lambda=\frac{\sqrt{1-R^2}}{R}$, and $c_1$, $c_2$ are arbitrary constants subject to $c_1^2-c_2^2=1$.

Thus, up to rigid motions in $\mathbb R^3$, we have a two-parametric family of circular tractrices: one parameter is the positive constant $R$ and another parameter is encoded in the constants $c_1$, $c_2$ related by one relation.
	
	\begin{figure}[ht]
		\begin{center}
			\includegraphics[scale=0.25]{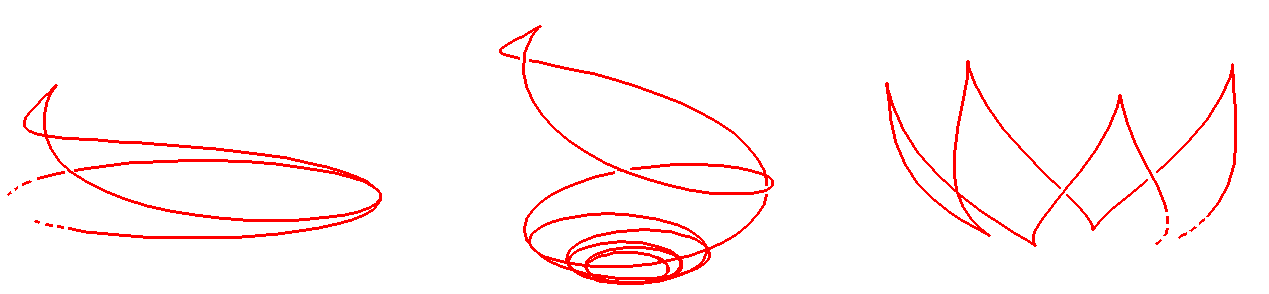}
			\caption{Typical examples of circular tractrices  with $R>1$, $R=1$, $R<1$}\label{fig1}
		\end{center}
	\end{figure}
	
	Qualitative properties of circular tractrices with $R$ greater, equal or less than 1 are strongly different. These three cases will be discussed separately in Chapters 2-4 below.


	\section{Circular tractrices with $R>1$}
	
	Fix $R>1$, chose arbitrary $c_1$, $c_2$ satisfying $c_1^2+c_2^2=1$, and consider the corresponding circular tractrix $\gamma$ represented by (\ref{pv}),(\ref{large}). Denote by $x = f(t)$ the position vector of $\gamma$. List elementary geometrical properties of $\gamma$.
	
	\smallskip
	
	1) The circular tractrix $\gamma$ is symmetric with respect to the plane $x^2=0$:
	\begin{eqnarray}\label{symmetry}
		f^1(-t) =  f^1(t),\quad f^2(-t) = - f^2(t),\quad f^3(-t) = f^3(t).
	\end{eqnarray}
	This is the unique symmetry of $\gamma$.
	
	\smallskip

	2) The circular tractrix  $\gamma$ is piecewise regular. It has a unique singular point, a cusp, at $t=0$. This follows immediately from the following relation:
	\begin{eqnarray}\label{2.1}
		\vert f^\prime \vert = \left\vert \frac{\lambda\sinh\lambda t}{\frac{c_1}{R}+\cosh\lambda t}\right\vert.
	\end{eqnarray}
	Notice that $\vert f^\prime \vert = \vert \xi_2 \vert $. Besides, by integrating (\ref{2.1}) one can show that $\gamma$ has infinite length.

	\smallskip
	
	3) As $t$ tends to $\pm\infty$, the circular tractrix $\gamma$ becomes asymptotically close to the circle $C_\infty$ of radius $\sqrt{R^2-1}$ centered at the origin in the coordinate plane $x^3=0$. More exactly, if we introduce two vector-functions,
	\begin{eqnarray*}
		f_+(t) =
		\lambda \left(
		\begin{array}{c}
			\sqrt{R^2-1}\cos\frac{t}{R} + \sin \frac{t}{R}\\
			\sqrt{R^2-1}\sin\frac{t}{R} - \cos \frac{t}{R}\\
			0
		\end{array}
		\right) , \quad
		f_- (t) =
		\lambda \left(
		\begin{array}{c}
			\sqrt{R^2-1}\cos\frac{t}{R} - \sin \frac{t}{R}\\
			\sqrt{R^2-1}\sin\frac{t}{R} + \cos \frac{t}{R}\\
			0
		\end{array}
		\right)
	\end{eqnarray*}
	which both represent the circle $C_\infty$ in appropriate parameterizations, then we have
	\begin{eqnarray*}
		\vert  f(t) - f_+ (t)\vert < 2 e^{-\lambda t}, \quad \vert  f(t) - f_- (t)\vert < 2 e^{\lambda t},
	\end{eqnarray*}
	and hence
	\begin{eqnarray*}
		\lim\limits_{t\to +\infty} \vert  f(t) - f_+ (t)\vert  = 0, \quad \lim\limits_{t\to -\infty} \vert  f(t) - f_- (t)\vert  = 0.
	\end{eqnarray*}
	Evidently, similar asymptotic closeness at $t\to\pm \infty$ extends to the Frenet frames, curvatures and torsions of $\gamma$ and $C_\infty$ respectively.
	
	Notice that the asymptotic circle $C_\infty$ does not depend on the choice of $c_1$, $c_2$.
	
	\smallskip
	
	4) The position vector $x=f(t)$ of the circular tractrix $\gamma$ satisfies the following relation:
	\begin{eqnarray}\label{tracing} f \, + \, \frac{1}{\xi_2} \, f^\prime =
		\left(
		\begin{array}{c}
			R\cos\frac{t}{R}\\ R\sin\frac{t}{R} \\ 0
		\end{array}
		\right).\end{eqnarray}
	
	This means that if one draws appropriately chosen unit segments tangent to $\gamma$, then the endpoints of these segments sweep out the circle $C$ of radius $R$ centered at the origin in the coordinate plane $x^3=0$, and  $t$  is an arc length for $C$. Thus, the circular tractrix $\gamma$ is related to the circle $C$ in the same manner as the classical {\it linear} tractix is related to its asymptotic straight line, c.f. \cite[p.8]{DoCarmo}. In terms of the general theory of tractrices, the circle $C$ is {\it the directrix} for the circular tractrix $\gamma$ in question, c.f. \cite{survey_tractrices}, \cite{survey_tractrices1}.
	Notice that the circle $C$ does not depend on the choice of $c_1$, $c_2$.
	
	\smallskip
	
	5) The circular tractrix $\gamma$ has non-vanishing torsion for any choice of $c_1$, $c_2$ except two particular cases, $c_1=1$, $c_2=0$ and $c_1= - 1$, $c_2=0$, where $\gamma$ belongs to the coordinate plane $x^3=0$ and represents the well known planar circular tractrices with $R>1$, see Fig.2, c.f. \cite{survey_tractrices}, \cite{survey_tractrices1}.
	
	\begin{figure}[ht]
		\begin{center}
			\includegraphics[scale=0.3]{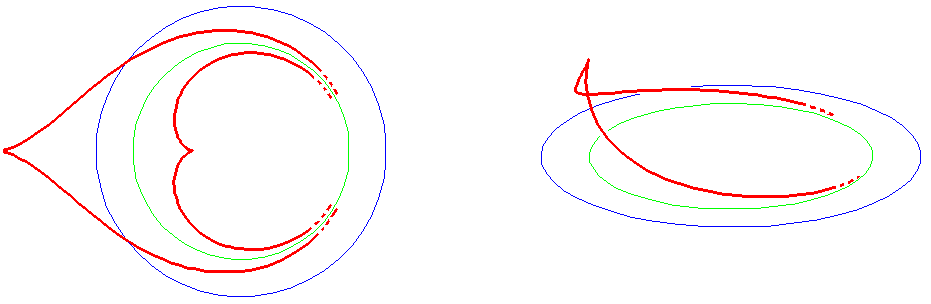}
			\caption{ Left: The planar circular tractrices (red), circle $C$ (blue) and circle $C_\infty$ (green) viewed from the top. Right: A non-planar circular tractrix (red) and the same circles $C$ (blue) and $C_\infty$ (green).}\label{fig2}
		\end{center}
	\end{figure}
	
	\smallskip
	
	6) The position vector $x=f(t)$ of the circular tractrix $\gamma$ satisfies the following relation:
	\begin{eqnarray}\label{pvlimit}
		\lim\limits_{R\to+\infty} \left(
		\begin{array}{c}
			f^1(t) - R \\ f^2(t)\\ f^3(t)
		\end{array}\right) = \left(
		\begin{array}{c}
			- c_1 \frac{1}{\cosh t},\\
			t - \tanh t,\\
			c_2 \frac{1}{\cosh t}.
		\end{array}\right) .
	\end{eqnarray}
	
	This means that if we use a shift along the $x^1$-axis at the distance $R$ so that the circle $C$ passes through the origin, then  at the limit $R\to+\infty$  the circle $C$ transforms into the $x^2$-axis, the circular tractrix $\gamma$ transforms into the well known linear tractrix situated so that its asymptotic straight line is the $x^2$-axis, and the parameters $c_1$, $c_2$ satisfying $c_1^2+c_2^2=1$ describe the rotation of the linear tractrix around the $x^2$-axis in $\mathbb R^3$, c.f. \cite[p.7]{Tenenblat}.
	
	\medskip

\begin{figure}[ht]
		\begin{center}
			\includegraphics[scale=0.3]{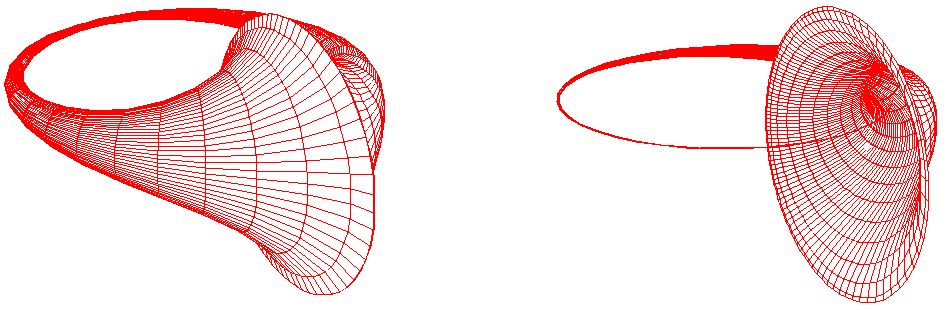}
			\caption{Typical example of circular pseudosphere with $R>1$: the complete surface (left) and its one half (right)}\label{fig3}
		\end{center}
	\end{figure}

	Next, the constant $R>1$ being fixed, set $c_1=\cos \alpha$, $c_2=\sin\alpha$ and allows $\alpha\in S^1$ to be varied. Then we obtain a one-parameter family of circular tractrices which sweep out a two-dimensional surface $F$. This surface is represented by the position vector $x=f(t,\alpha)$ given by (\ref{pv}), (\ref{large}) with $c_1=\cos \alpha$, $c_2=\sin\alpha$. We will call $F$ {\it a circular pseudosphere}, see Fig.\ref{fig3}.

	Let us list fundamental geometric properties of $F$.
	
	\smallskip
	
	1*) The circular pseudosphere $F$ is symmetric with respect to the coordinate plane $x^2=0$:
	\begin{eqnarray}\label{symmetryF}
		f^1(-t,\alpha) =  f^1(t,\alpha),\quad f^2(-t,\alpha) = - f^2(t,\alpha),\quad f^3(-t,\alpha) = f^3(t,\alpha).
	\end{eqnarray}
	Thus, $F$ consists of 2 symmetric parts sharing the common coordinate line $t=0$ situated in the plane $x^2=0$.
	
	Moreover, $F$ is symmetric with respect to the coordinate plane $x^3=0$:
	\begin{eqnarray}\label{symmetryFF}
		f^1(t,-\alpha) =  f^1(t,\alpha),\quad f^2(t,-\alpha) = f^2(t,\alpha),\quad f^3(t,-\alpha) = - f^3(t,\alpha).
	\end{eqnarray}
	
	\smallskip
	
	2*) The circular pseudosphere $F$ is piecewise regular. Its singular set, a cuspidal edge, is a unit circle composed of the singular points $t=0$ of the circular tractrices $\alpha=const$ sweeping out the surface $F$.
	
	\smallskip
	
	3*) The coordinate curves $t=const$ in $F$ are circles whose radii are equal to $\frac{\sqrt{R^2-1}}{\sqrt{R^2\cosh^2\lambda t -1}}$ and tend to $0$ as $t\to\pm\infty$. This is verified easily by computing the curvature and torsion of the curves in question viewed as curves in $\mathbb R^3$.
	
	\smallskip
	
	4*) As $t$ tends to $\pm\infty$, the circular pseudosphere $F$  becomes asymptotically close to the circle $C_\infty$ in the same manner as it was described above for the circular tractrices constituting $F$.
	
	\smallskip

	5*) The position vector $x=f(t,\alpha)$ of the circular pseudosphere $F$ satisfies the following relation:
	\begin{eqnarray}\label{position_vector_s}
		f \, + \, \frac{1}{\xi_2} \cdot \frac{\partial f}{\partial t} =
		\left(
		\begin{array}{c}
			R\cos\frac{t}{R}\\ R\sin\frac{t}{R} \\ 0
		\end{array}
		\right).
	\end{eqnarray}
	
	This means that if one draws appropriately chosen unit segments tangent to coordinate $t$-curves of $F$, then the endpoints of these segments will sweep out the circle $C$. Thus, the circular pseudosphere $F$ is related to the circle $C$ in the same manner as the classical pseudosphere is related to its axis of rotation.
	
	\smallskip
	
	6*) The first fundamental form of $F$ reads as follows:
	\begin{eqnarray*}
		g = \frac{R^2-1}{(\cos\alpha + R\cosh\lambda t)^2} \left( \sinh^2\lambda t \, dt^2 \, + \, d\alpha^2 \right).
	\end{eqnarray*}
	Hence, the coordinate curves in $F$, which are the circular tractrices $\alpha=const$ and the circles $t=const$, form an isothermic net on $F$. Evidently, this isothermic net on the circular tractrix $F$ can be viewed as an analogue of the standard horocyclic net on the classical pseudosphere.
	
	\smallskip
	
\begin{figure}[ht]
		\begin{center}
			\includegraphics[scale=0.4]{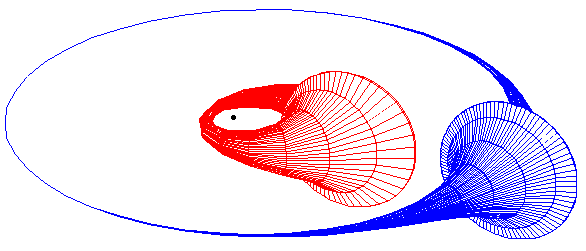}
			\caption{Circular pseudospheres with $R=1.25$ (red) and $R=3.5$ (blue)}\label{fig4}
		\end{center}
	\end{figure}

	7*) Clearly, the circular pseudosphere $F$ depends on $R$. For instance, it is  situated inside the tube of unit radius around the circle  $C$ in $\mathbb R^3$.  Hence, the greater $R$ is, the more distant $F$ is from the origin $O\in\mathbb R^3$, see Fig.\ref{fig4}.
	
	\smallskip
	
	8*) No matter what $R>1$ is, the complete area of the circular pseudosphere $F$ is equal to $4\pi$, and it is the same as the area of the classical pseudosphere, c.f. \cite{psbook}. This is a really astonishing fact, in view of the previous item in the list. Its proof is based simply on calculations of corresponding integrals:
	\begin{eqnarray*}
		Area(F) = 2\int\limits_{0}^{+\infty} \int\limits_{0}^{2\pi} \sqrt{det g}\, d\alpha\, dt \, = \, 2 \int\limits_{0}^{+\infty} \int\limits_{0}^{2\pi} \frac{(R^2-1)\sinh\lambda t}{(\cos\alpha + R\cosh\lambda t)^2}\, d\alpha\, dt \, =\, 4\pi.
	\end{eqnarray*}
	
	\smallskip
	
	9*) The circular pseudosphere $F$ has self-intersections, which all are situated in the coordinate plane $x^2=0$. By cutting $F$ with this plane, we  decompose $F$ into a countable set of pieces without self-intersections, $F=\cup_{j=1}^\infty F_j$. Every piece $F_j$ encloses a well-defined body, $F_j=\partial B_j$. Then we define "the volume enclosed by" $F$ as the total sum of the volumes of $B_j$, $j\geq 1$, i.e., $Vol(F):= \sum\limits_{j=1}^\infty Vol(B_j)$. Actually, $Vol(F)$ is the volume of a body enclosed by $F$, which is counted with multiplicities  in view of self-intersections of $F$.
	
	To find $Vol(F)$, one needs to parameterize every body $B_j$. For instance, $F_j$ is equipped with the coordinates $t$, $\alpha$ and foliated by circles $t=const$. Consequently, $B_j$ is foliated by two-dimensional discs. Introducing polar coordinates $l$, $\alpha$ in every disc, we get a parameterization by $t$, $\alpha$ $l$ for $B_j$.  And then the volume of $B_j$ is calculated via appropriate integrals in terms of $t$, $\alpha$ $l$, which appear to be quite cumbersome. However the result turns out to be absolutely surprising:
	\begin{eqnarray}\label{volume}
		Vol(F) \, = \frac{2}{3}\pi,
	\end{eqnarray}
	Thus, no matter what $R>1$ is, the complete "volume enclosed by"\, the circular pseudosphere $F$ is equal to $\frac{2}{3}\pi$, and it is the same as the volume enclosed by the classical pseudosphere, c.f. \cite{psbook}.
	
	\smallskip
	
	10*) If we make a shift along the $x^1$-axis at the distance $R$ so that the circle $C$ passes through the origin, then  at the limit $R\to +\infty$  the circle $C$ transforms into the $x^2$-axis, the circular pseudosphere $F$ transforms into the classical pseudosphere whose axis of rotation is the $x^2$-axis, and $\alpha$ describes the corresponding angle parameter on the pseudosphere, c.f. \cite[p.7]{Tenenblat}, \cite{psbook}. Therefore, the pseudosphere appears as the limit surface at $R\to \infty$ in the one-parametric family of circular pseudospheres with $R>1$ under consideration.

	\medskip
	
	Thus, the circular tractrices and circular pseudospheres with $R>1$ inherit fundamental properties of the classical linear tractrix and pseudosphere respectively and hence can be naturally treated as their circular analogs.

	\section{Circular tractrices with $R=1$}

	Next fix $R=1$, chose arbitrary $c_1$, $c_2$ satisfying $4(c_1-1) = c_2^2$, and consider the corresponding circular tractrix $\gamma$ represented by (\ref{pv}),(\ref{unit}). Let $x = f(t)$ stand for the position vector of $\gamma$. List elementary geometrical properties of $\gamma$.
	
	\smallskip
	
	1) $\gamma$ is symmetric with respect to the plane $x^2=0$, its position-vector satisfies (\ref{symmetry}).
	This is the unique symmetry of $\gamma$.
	
	\smallskip

	2) $\gamma$ is piecewise regular. It has a unique singular point, a cusp, at $t=0$. This follows immediately from the relation:
	\begin{eqnarray*}
		\vert f^\prime \vert = \left\vert \frac{2t}{t^2+c_1} \right\vert.
	\end{eqnarray*}
	As consequence,  $\gamma$ has infinite length. Notice that $\vert f^\prime \vert = \vert \xi_2 \vert $ remains true.

	\smallskip
	
	3) As $t$ tends to $\pm\infty$, the circular tractrix $\gamma$ becomes asymptotically close to the origin point $O\in\mathbb R^3$. More exactly, we have
	\begin{eqnarray*}
		\vert  f(t) \vert < \frac{2}{\vert t\vert},
	\end{eqnarray*}
	and hence $\vert  f(t) \vert  \to 0$ as $t\to\pm\infty$. The origin point $O$ can be viewed as a degenerate asymptotic circle $C_\infty$ whose radius $\sqrt{R^2-1}$ become equal to $0$ at $R=1$. Notice that this asymptotic behavior  of $\gamma$ does not depend on the choice of $c_1$, $c_2$.
	
	\smallskip
	
	4) The position vector $x=f(t)$ of the circular tractrix $\gamma$ satisfies the same relation as (\ref{tracing}) but with $R=1$. This still means that if one draws appropriately chosen unit segments tangent to $\gamma$, then the endpoints of these segments will sweep out the circle $C$ of unit radius centered at the origin in the coordinate plane $x^3=0$, and  $t$  is an arc length of this circle. Thus, the circular tractrix $\gamma$ is related to the unit circle $C$ in the same manner as the linear tractix is related to its asymptotic straight line, i.e., the unit circle $C$ is the directrix for the circular tractrix $\gamma$ in question, c.f. \cite[p.8]{DoCarmo}, \cite{survey_tractrices}, \cite{survey_tractrices1}. Notice that the circle $C$ does not depend on the choice of $c_1$, $c_2$.
	
	\smallskip
	
	5) The circular tractrix $\gamma$ has non-vanishing torsion for any choice of $c_1$, $c_2$ except one particular case, $c_1=1$, $c_2=0$, where $\gamma$ belongs to the coordinate plane $x^3=0$ and describes the well known planar circular tractrix with $R=1$, see Fig.5, c.f. \cite{survey_tractrices}, \cite{survey_tractrices1}.
	
	\begin{figure}[ht]
		\begin{center}
			\includegraphics[scale=0.3]{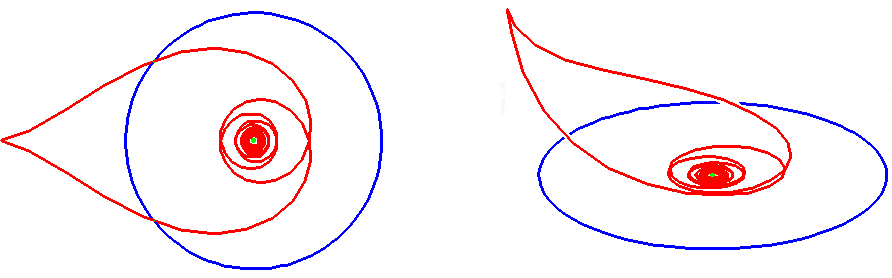}
			\caption{ Left: The planar circular tractrix (red) and the unit circle $C$ (blue) viewed from the top. Right: A non-planar circular tractrix (red) and the unit circle $C$ (blue).   }\label{fig5}
		\end{center}
	\end{figure}
	
\begin{figure}[ht]
		\begin{center}
			\includegraphics[scale=0.25]{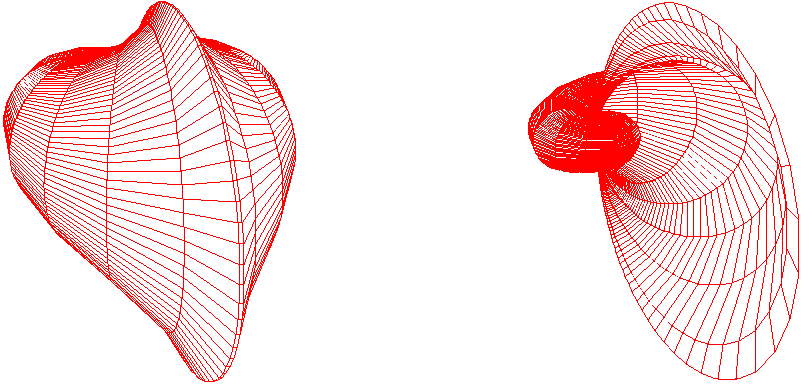}
			\caption{The circular pseudosphere with $R=1$: complete surface (left) and its one half (right)}\label{fig6}
		\end{center}
	\end{figure}

	Next, the constant $R=1$ being fixed, set $c_1=1+\alpha^2$, $c_2=2\alpha$ and allows $\alpha\in \mathbb R^1$ to be varied. Then we obtain a one-parameter family of circular tractrices, which sweep out a two-dimensional surface $F$. This surfaces, which is also called {\it a circular pseudosphere}, is represented by the position vector $x=f(t,\alpha)$ given by (\ref{pv}),(\ref{unit}) with $c_1=1+\alpha^2$, $c_2=2\alpha$ , see Fig.\ref{fig6}.
	
	Let us list fundamental geometric properties of $F$.
	
	\smallskip
	
	1*) The circular pseudosphere $F$ is symmetric with respect to the coordinate plane $x^2=0$, its position vector satisfies (\ref{symmetryF}). Thus, $F$ consists of two symmetric parts sharing the common coordinate line $t=0$ situated in the plane $x^2=0$. Besides, $F$ is symmetric with respect to the coordinate plane $x^3=0$, its position vector satisfies (\ref{symmetryFF}).
	
	\smallskip
	
	2*) The circular pseudosphere $F$ is piecewise regular. Its singular set, a cuspidal edge, is a unit circle composed of the singular points $t=0$ of the circular tractrices sweeping out the surface $F$.
	
	\smallskip
	
	3*) The coordinate curves $t=const$ in $F$ are circles of radius $\frac{1}{\sqrt{t^2+1}}$. All of them pass through the origin point $O\in \mathbb R^3$ which corresponds to the limit value $\alpha \to \pm\infty$ and has to be treated as removed from $F$.
	
	\smallskip
	
	4*) As $t$ tends to $\pm\infty$, the circular pseudosphere $F$  becomes asymptotically close to the origin $O\in \mathbb R^3$.
	
	\smallskip
	
	5*) The position vector $x=f(t,\alpha)$ of the circular pseudosphere $F$ satisfies the same relation as in (\ref{position_vector_s}) but with $R=1$. Hence, if one draws appropriately chosen unit segments tangent to coordinate $t$-curves of $F$, then the endpoints of these segments will sweep out the unit circle $C$. Moreover, $t$ is still an arc length of this circle. Thus, the circular pseudosphere $F$ in question is related to the unit circle $C$ in the same manner as the classical pseudosphere is related to its axis of rotation.
	
	\smallskip
	
	6*) The first fundamental form of $F$ reads as follows
	\begin{eqnarray*}
		g = \frac{4}{(1+\alpha^2 + t^2)^2} \left( t^2 \, dt^2 \, + \, d\alpha^2 \right).
	\end{eqnarray*}
	Hence, the coordinate curves in $F$, circular tractrices $\alpha=const$ and circles $t=const$, form an isothermic net on $F$. Possibly, this isothermic net on the circular tractrix $F$ may be viewed as analogue of the standard horocyclic net on the classical pseudosphere.
	
	\smallskip
	
	7*) The complete area of the circular pseudosphere $F$ is still equal to $4\pi$, and it is the same as the area of the classical pseudosphere. The proof is based on calculations of corresponding integrals:
	$$
	Area(F) = 2\int\limits_{0}^{+\infty} \int\limits_{-\infty}^{+\infty} \sqrt{det g}\, d\alpha\, dt \, = \, 2 \int\limits_{0}^{+\infty} \int\limits_{-\infty}^{+\infty} \frac{4 t}{(1+\alpha^2 + t^2)^2}\, d\alpha\, dt \, =\, 4\pi.
	$$
	
	\smallskip
	
	8*) The circular pseudosphere $F$ has self-intersections. We can define the "volume enclosed by"\, $F$  in the same manner as we use in the previous section for the case of $R>1$. Then we get the same formula (\ref{volume}).	 Therefore, the complete "volume enclosed by"\, the circular pseudosphere $F$ is still equal to $\frac{2}{3}\pi$, and it is the same as the volume enclosed by the classical pseudosphere.
	
	\medskip
	
	Thus, similarly to the case  of $R>1$, the circular tractrices and circular pseudosphere with $R=1$ inherit fundamental properties of the classical linear tractrix and pseudosphere respectively and hence can be treated as their circular analogs too.
	
	Notice that the case $R=1$ discussed in this section can be viewed as the limit for the case $R>1$ considered in the previous section. Particularly, formulae (\ref{unit}) arise as the limit version of (\ref{large}) as $R\to 1$ under appropriate scalings of involved parameters. Besides, qualitative geometric properties of circular tractrices and circular pseudospheres with $R>1$ hold in the limit case $R=1$.

	\section{ Circular tractrices with $0<R<1$ }

	Finally fix $0<R<1$, chose arbitrary $c_1$, $c_2$ satisfying $c_1^2 - c_2^2=1$, and consider the corresponding circular tractrix $\gamma$ represented by (\ref{pv}),(\ref{small}). Let $x = f(t)$ still stand for the position vector of $\gamma$. List elementary geometrical properties of $\gamma$.
	
	\smallskip
	
	1) $\gamma$ is symmetric with respect to the plane $x^2=0$, its position vector satisfies (\ref{symmetry}).
	Moreover, $\gamma$ is periodic in the following sense:
	\begin{eqnarray*}
		f^1(t+T) =  f^1(t) \cos\varphi + f^2(t) \sin\varphi,
	\end{eqnarray*}
	\begin{eqnarray*}
		f^2(t+T) = - f^1(t) \sin\varphi + f^2(t) \cos\varphi , \quad f^3(t+T) = f^3(t),
	\end{eqnarray*}
	where $T=\frac{2\pi}{\lambda}$ and $\varphi=\frac{2\pi }{\sqrt{1-R^2}}$. Thus, $\gamma$ is invariant under rotations around the $x^3$-axis at the angles $n\varphi$, $n\in\mathbb Z$, in $\mathbb R^3$. Consequently, $\gamma$ is symmetric with respect to any plane $x^1 \sin\frac{\varphi}{2}n + x^2 \cos\frac{\varphi}{2}n = 0$, $n\in\mathbb Z$, in $\mathbb R^3$.
	
	\smallskip

	2) $\gamma$ is piecewise regular. It has two rotationally invariant series of singular points: $t=\frac{2n\pi}{\lambda}$, $n\in\mathbb Z$, and $t=\frac{(2n+1)\pi}{\lambda}$, $n\in\mathbb Z$. This follows immediately from the relation:
	\begin{eqnarray*}
		\vert f^\prime \vert = \left\vert \frac{\lambda\sin\lambda t}{\frac{c_1}{R}+\cos\lambda t} \right\vert.
	\end{eqnarray*}
	Notice that $\vert f^\prime \vert = \vert \xi_2 \vert $ remains true.
	
	\smallskip
	
\begin{figure}[ht]
		\begin{center}
			\includegraphics[scale=0.25]{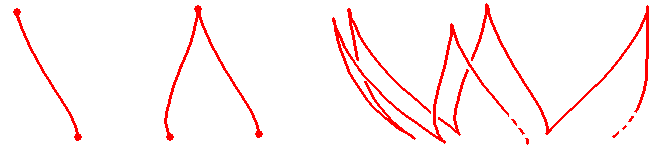}
			\caption{A circular tractrix with $R<1$ (right), its  unit (left) and petal (center)}\label{fig7}
		\end{center}
	\end{figure}

	3) Any interval of $\gamma$, which is situated between two consecutive singular points, will be called {\it a unit} of $\gamma$. Any pair of adjacent units represents a piece of $\gamma$, which is symmetric with respect to the two-dimensional plane containing the $x^3$-axis and the common endpoint of both units; it will be called {\it a petal} of $\gamma$, see Fig.\ref{fig7}. The complete circular tractrix  $\gamma$ is reconstructed by applying to any of its petals the rotations around the $x^3$-axis at the angles $n\varphi$, $n\in\mathbb Z$.
	
	\smallskip
	
	Particularly, $\gamma$ is closed if and only if $\frac{\varphi}{\pi}\in \mathbb Q$, i.e. $\sqrt{1-R^2}\in\mathbb Q$. Otherwise, $\gamma$ forms an everywhere dense subset in some rotationally invariant surface in $\mathbb R^3$.
	
	\begin{figure}[ht]
		\begin{center}
			\includegraphics[scale=0.25]{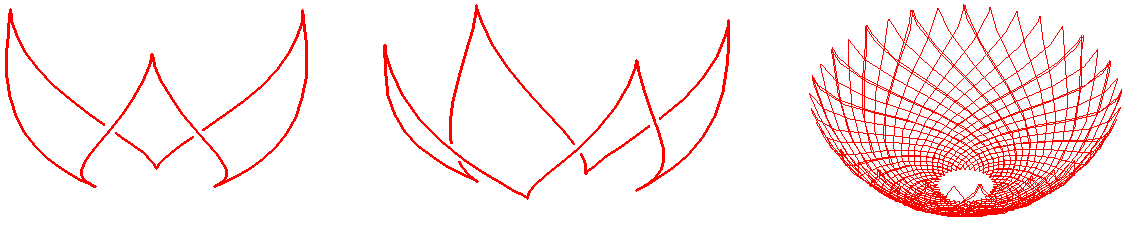}
			\caption{Closed circular tractrix with $R=\sqrt{1-\left(\frac{3}{4}\right)^2}$ (left) and $R=\sqrt{1-\left(\frac{4}{5}\right)^2}$ (center), and a non-closed circular tractrix with $R=\sqrt{1-\frac{1}{e}}$ (right).}\label{fig8}
		\end{center}
	\end{figure}
	
	\smallskip
	
	4) The length of a unit of $\gamma$ is equal to $sign(c_1) \ln\left\vert\frac{c_1+R}{c_1-R}\right\vert$. Hence it depends on $R$ as well as on $c_1$. Notice that the length is finite, no matter what $0<R<1$ and $c_1\in \mathbb R$ are. But it can tend to infinity as $R\to 1$, and this illustrates a quite complex behaviour of the circular tractrix $\gamma$ as $R$ approaches 1 from below.
	
	\smallskip
	
	5) The position vector $x=f(t)$ of the circular tractrix $\gamma$ satisfies the same relation as in (\ref{tracing}) but with $0<R<1$. Once again, this means that if one draws appropriately chosen unit segments tangent to $\gamma$, then the endpoints of these segments will sweep out the circle $C$ of radius $R$ centered at the origin in the coordinate plane $x^3=0$, and  $t$  is an arc length of this circle. Thus, the circular tractrix $\gamma$ is related to the circle $C$ in the same manner as the linear tractix is related to its asymptotic straight line, i.e. the circle $C$ is the directrix for the circular tractrix $\gamma$ in question, c.f. \cite[p.8]{DoCarmo}, \cite{survey_tractrices}, \cite{survey_tractrices1}. Notice that the circle $C$ does not depend on the choice of $c_1$, $c_2$.
	
	\smallskip
	
	6) The circular tractrix $\gamma$ has non-vanishing torsion for any choice of $c_1$, $c_2$ except two particular cases, $c_1=1$, $c_2=0$ and $c_1= - 1$, $c_2=0$, where $\gamma$ belongs to the coordinate plane $x^3=0$ and represents the well known planar circular tractrices with $0<R<1$, see Fig.9, c.f. \cite{survey_tractrices}, \cite{survey_tractrices1}.
	
	\begin{figure}[ht]
		\begin{center}
			\includegraphics[scale=0.25]{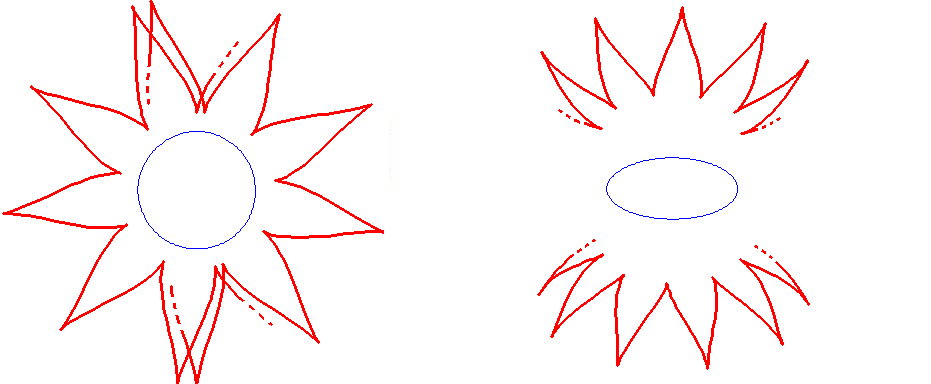}
			\caption{ Left: The planar circular tractrices (red) and circle $C$ (blue) viewed from the top. Right: Non-planar circular tractrices (red) and the circle $C$ (blue)}\label{fig9}
		\end{center}
	\end{figure}

	\medskip

	Next, the constant $0<R<1$ being fixed, set $c_1=\pm\cosh \alpha$, $c_2=\sinh\alpha$ and allows $\alpha\in \mathbb R^1$ to be varied. Then we obtain a one-parameter two-component family of circular tractrices which sweep out a two-dimensional surface $F$. This surface is represented by the position vector $x=f(t,\alpha)$ given by (\ref{pv}),(\ref{small}) with $c_1=\pm\cosh \alpha$, $c_2=\sinh\alpha$, it will be called {\it a circular pseudosphere} too. Notice that the surface $F$ consists of two mutually congruent components which correspond to  $c_1 = \cosh \alpha$ and $c_1 = - \cosh \alpha$ respectively.

	Let us list fundamental geometric properties of $F$.
	
	\smallskip
	
	1*) The circular pseudosphere $F$ is symmetric with respect to the coordinate plane $x^2=0$, its position vector satisfies (\ref{symmetryF}). Moreover,  $F$ is periodic in the following sense:
	\begin{eqnarray*}
		f^1(t+T,\alpha) =  f^1(t,\alpha) \cos\varphi + f^2(t,\alpha) \sin\varphi,
	\end{eqnarray*}
	\begin{eqnarray*}
		f^2(t+T,\alpha) = - f^1(t,\alpha) \sin\varphi + f^2(t,\alpha) \cos\varphi , \quad f^3(t+T,\alpha) = f^3(t,\alpha),
	\end{eqnarray*}
	where $T=\frac{2\pi}{\lambda}$ and $\varphi=\frac{2\pi }{\sqrt{1-R^2}}$. Thus, the surface $F$ is invariant under rotations around the $x^3$-axis at the angles $n\varphi$, $n\in\mathbb Z$, in $\mathbb R^3$. Consequently, $F$ is symmetric with respect to any plane $x^1 \sin\frac{\varphi}{2}n + x^2 \cos\frac{\varphi}{2}n = 0$, $n\in\mathbb Z$, in $\mathbb R^3$. Besides, $F$ is symmetric with respect to the coordinate plane $x^3=0$, its position vector satisfies (\ref{symmetryFF}).
	
	\smallskip

	2*) The circular pseudosphere $F$ is piecewise regular. It has two rotationally invariant series of  cuspidal edges, the coordinate curves $t=\frac{2n\pi}{\lambda}$, $n\in\mathbb Z$, and $t=\frac{(2n+1)\pi}{\lambda}$, $n\in\mathbb Z$, respectively, which are formed by singular points of circular tractrices $\alpha=const$ sweeping out the surface $F$.
	
	\smallskip
	
	3*) The coordinate curves $t=const$ in $F$ are circles of radii $\frac{\sqrt{1-R^2}}{\sqrt{1-R^2\cos^2\lambda t}}$, this fact is easily verified by computing the curvature and torsion of the mentioned curves viewed as curves in $\mathbb R^3$. Particularly, singular edges of $F$ are unit circles.
	
	\smallskip
	
	4*) All the coordinate circles $t=const$ of $F$ pass through the points $O_1(0,0,-\sqrt{1-R^2})$ and $O_2(0,0,\sqrt{1-R^2})$, which correspond to the limit values $\alpha\to\pm \infty$ and hence have to be viewed as removed from $F$.
	
	\smallskip
	
\begin{figure}[ht]
		\begin{center}
			\includegraphics[scale=0.25]{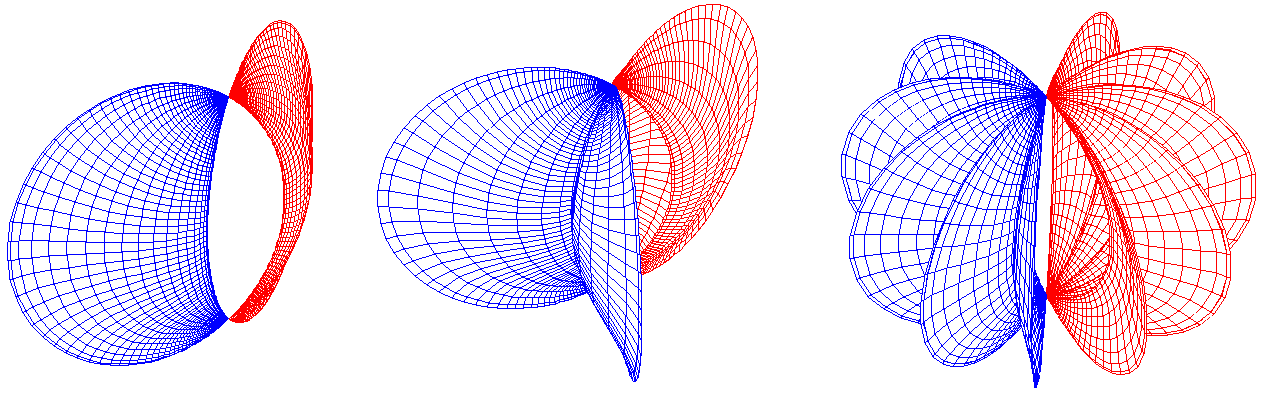}
			\caption{A circular pseudosphere with $R<1$:  a unit (left), a petal (center), several petals (right). Different components are colored by red and blue. }\label{fig10}
		\end{center}
	\end{figure}

	5*) Any piece of $F$, which is determined by $\frac{n\pi}{\lambda} \leq t\leq \frac{(n+1)\pi}{\lambda}$ for some $n\in\mathbb Z$ and hence situated between two consecutive singular edges, will be called {\it a unit} of $F$. Any pair of adjacent units represents a piece of $F$, which is symmetric with respect to the two-dimensional plane containing the $x^3$-axis and the singular coordinate circle shared by the units in question. This piece will be called {\it a petal} of $F$, see Fig.\ref{fig10}. The complete circular pseudosphere  $F$ is reconstructed by applying to any of its petals the rotations around the $x^3$-axis at the angles $n\varphi$, $n\in\mathbb Z$.
	
	\smallskip
	
	6*) Any unit of $F$ has two connected components, which corresponds to the choice either $c_1=\cosh\alpha$ or $c_1=-\cosh\alpha$. Both components share the same pair of points  $O_1(0,0,-\sqrt{1-R^2})$ and $O_2(0,0,\sqrt{1-R^2})$ viewed as removed from $F$. Moreover, these components are mutually congruent.

	\smallskip
	
	7*) The position vector $x=f(t,\alpha)$ of the circular pseudosphere $F$ satisfies the same relation as in (\ref{position_vector_s}) but with $0<R<1$. This still means that if one draws appropriately chosen unit segments tangent to coordinate $t$-curves of $F$, then the endpoints of these segments will sweep out the circle $C$. Moreover, $t$ is an arc length of this circle. Thus, the circular pseudosphere $F$ is related to the circle $C$ in the same manner as the classical pseudosphere is related to its axis of rotation.

	\smallskip
	
	8*) The first fundamental form of $F$ reads as follows:
	\begin{eqnarray*}
		g = \frac{1-R^2}{(\pm\cosh\alpha + R\cos\lambda t)^2} \left( \sin^2\lambda t \, dt^2 \, + \, d\alpha^2 \right).
	\end{eqnarray*}
	Hence, the coordinate curves in $F$, which are the circular tractrices $\alpha=const$ and the circles $t=const$, form an isothermic net on $F$ which can be viewed as an analogue of the standard horocyclic net on the classical pseudosphere.
	
	\smallskip
	
	9*) The complete area of a unit of $F$ is equal to
	\begin{eqnarray*}
		Area = 2\int\limits_{-\infty}^{+\infty} \int\limits_{0}^{\frac{\pi}{\lambda}} \sqrt{det g}\, dt \, d\alpha\,  = \, 2\int\limits_{-\infty}^{+\infty} \int\limits_{0}^{\frac{\pi}{\lambda}} \frac{(1-R^2)\sin\lambda t}{(\cosh\alpha + R\cos\lambda t)^2}\, dt \, d\alpha\,  = \\ = \, 4\left( \arctan\sqrt{\frac{1+R}{1-R}}-\arctan\sqrt{\frac{1-R}{1+R}} \right).
	\end{eqnarray*}
	Therefore, the area is finite and depends on $R$. In this context the case of $0<R<1$ differs essentially from the case of $R\geq 1$. Similar differences appear as well in the context of the volume enclosed by the circular pseudosphere.

	\section{Concluding remarks and questions}
	
	{\it Remark 1.} Notice that the Gauss curvature of circular pseudospheres is not constant negative, no matter what $R>0$ is. In other words, the circular pseudospheres are not pseudospherical in the classical sense of this term. On the other hand, fix $R>0$ and consider the corresponding one-parametric family of circular tractrices. All these tractrices have the same directrix, the circle $C$. By applying rotations along $C$ in $\mathbb R^3$, we get a two-parametric family of circular tractrices with the same directrix $C$. The question is whether one can choose a one-parametric subfamily in this two-parametric family of circular tractrices so that the chosen circular tractrices sweep out a surface of constant negative curvature. If such a pseudospherical surface exists, then it can be treated as a circular analog of the well-known Dini surface.
	
	{\it Remark 2.} It would be interesting to explore the extrinsic geometry of the circular pseudospheres. For instance, we claim that if one considers an arbitrary circular pseudosphere $F\subset\mathbb R^3$ parameterized by coordinates ($t$,$\alpha$) we use above, then the coordinate lines are {\it lines of curvature} in $F$. Moreover, the asymptotic lines in $F$ which live near singular edges of $F$ turn out to traverse these edges tangentially. In this regard circular pseudospheres show resemblance with the classical pseudosphere too. Possibly, another surprising resemblances can be found in this context.
	
	{\it Remark 3.} An arbitrary non-closed circular tractrix with $0<R<1$ form an everywhere dense subset in a rotationally invariant surface in $\mathbb R^3$. What can we say about that surface of revolution?
	
	{\it Remark 4.}  The circular tractrices can be characterized, with some exceptions, as the only tractrices in $\mathbb R^3$ whose directrices are circles. Namely, if we fix a circle $C\subset \mathbb R^3$ of radius $R$ then the tractrices in $\mathbb R^3$ whose directrix is $C$ are:
	
	$R>1$) the circular tractrices whose directix is $C$ and, as an exception, their common asymptotic circle $C_\infty$;
	
	$R=1$) the circular tractrices whose directix is $C$ and, as a degenerate exception, the origin point $O\in R^3$;
	
	$0<R<1$) the circular tractrices whose directix is $C$ and, as a degenerate exception, the corresponding points $O_1$, $O_2\in R^3$ that are situated from the {\it all} points of $C$ at the same distance 1.
	
	{\it Remark 5.} An interesting problem, which seems to be quite non-trivial, is to find helical analogues  for linear and circular tractrices in $\mathbb R^n$, $n\geq 3$. Namely, describe explicitly the tractrices in $\mathbb R^{n}$ whose directrices are curves of constant curvatures.

	{\it Remark 6.} In terms of the simple model of bicycle motion discussed in \cite{tabachnikov}, \cite{tabachnikov1}, any circular tractix in $\mathbb R^3$ can be viewed as the rear track of a spatial bicycle of unit length whose front track is a circle. In this context, the explicit description of the circular tractrices explored in our research note could be used for illustrating deep mathematical ideas and statements from \cite{tabachnikov}, \cite{tabachnikov1} concerning  tractrices in the frames of the modern theory of integrable systems.



\end{document}